\documentclass[11pt]{amsart}

\usepackage{amsmath,amssymb,amscd,latexsym}
\usepackage[sans]{dsfont}

\numberwithin{equation}{section}	% add the section number to the equation label

\title{Algebraic families of subfields \\ in division rings}

\author[J.-M. Bois, \lowercase{and} G. Vernik]{Jean-Marie Bois \lowercase{and} Gil Vernik}

\address{Mathematisches Seminar, Christian-Albrechts-Universit\"at zu Kiel, Ludewig-Meyn-Str.~4, 24098 Kiel, Germany}
\email{bois@math.uni-kiel.de, vernik@math.uni-kiel.de}
\thanks{The first named author is supported by the D.F.G. priority programme SPP1388 ``Darstellungstheorie''. The second named author is supported by Minerva Fellowship Program.}

\date{\today}


\newenvironment{thm}{\noindent \bf Theorem.\it}{ \rm }
\newenvironment{prop}{\noindent \bf Proposition.\it}{ \rm }
\newenvironment{lemma}{\noindent \bf Lemma.\it}{ \rm }
\newenvironment{cor}{\noindent \bf Corollary.\it}{ \rm }

\renewcommand{\proof}{\noindent {\it Proof.} }

\newcommand{\Frac}{\mathrm{Frac}}
\newcommand{\id}{\mathrm{id}}
\newcommand{\Char}{\mathrm{char}}
\newcommand{\ad}{\mathrm{ad}}
\newcommand{\Der}{\mathrm{Der}}
\newcommand{\Hom}{\mathrm{Hom}}

\newcommand{\Gal}{\mathrm{Gal}}
\newcommand{\Aut}{\mathrm{Aut}}
\newcommand{\Exp}{\mathrm{Exp}}
\newcommand{\Br}{\mathrm{Br}}

\newcommand{\cB}{\mathcal{B}}

\newcommand{\FF}{\mathbb{F}}
\newcommand{\ZZ}{\mathbb{Z}}

\begin{document}

\renewcommand{\labelitemi}{\textbullet}

\makeatletter \def\revddots{\mathinner{\mkern1mu\raise\p@ \vbox{\kern7\p@\hbox{.}}\mkern2mu \raise4\p@\hbox{.}\mkern2mu\raise7\p@\hbox{.}\mkern1mu}} \makeatother

\begin{abstract}
We describe relations between maximal subfields in a division ring and in its rational extensions. More precisely, we prove that properties such as being Galois or purely inseparable over the centre generically carry over from one to another. We provide an application to enveloping skewfields in positive characteristics. Namely, there always exist two maximal subfields of the enveloping skewfield of a solvable Lie algebra, such that one is Galois and the second purely inseparable of exponent 1 over the centre. This extends results of Schue in the restricted case \cite{SchuePaper}. Along the way we provide a description of the enveloping algebra of the $p$-envelope of a Lie algebra as a polynomial extension of the smaller enveloping algebra.
\end{abstract}

\maketitle

\section*{Introduction}
Let $D$ be a division ring which is finitely generated over its centre $Z$, and let $K$ be a subfield of $D$. The centralizer of $K$ in $D$ is defined by $C_{D}(K) = \{x \in D\ |\ [x,K] = 0\}$. The subfield $K$ is called a {\it maximal subfield of $D$} if $C_{D}(K) = K$. Alternatively, a subfield $K \subseteq D$ containing $Z$ is maximal if and only if  $[D:K] = [K:Z] = \sqrt{[D:Z]}$ \cite[Thm. 4.2.2 and 4.3.2]{hbook}. For more details about maximal subfields in the division rings one is referred to \cite{hbook,rowenBook}.  \\

A natural question is whether any central simple algebra affords a maximal subfield which is Galois over the centre (equivalently, whether such an algebra is a crossed product). The answer to this question is negative in general, see \cite[Thm. 7.1.30]{rowenBook}. In some special cases the answer may be positive. In \cite{SchuePaper}, J. Schue showed that this is the case for the division ring of fractions of the enveloping algebra of a solvable Lie $p$-algebra over a field $\FF$ of characteristic $p > 2$. In addition, in \cite{SchuePaper} was also shown the existence of a maximal subfield which is purely inseparable of exponent one over the centre.  \\

The present paper is concerned with the compared structures of maximal subfields in a division ring $D$ and in the division ring of rational functions $D(X)$ (see Section \ref{rationalfunctions} for a formal definition). Maximal subfields of $D(X)$ can be interpreted as parametrised families of subfields in $D$. The main result of this paper is a transfer theorem showing in what way properties of an extension $K:Z$ (such as being Galois or purely inseparable) pass from the corresponding properties for maximal subfields of $D(X)$ (Theorem \ref{reduction}). \\

The paper is organised as follows. In Section \ref{part1}, we assume that $D$ is a $p$-division algebra, ie.\ the dimension $[D:Z]$ is a power of $p = \Char(Z)$. In that situation, we provide a link between the notions of tori in $D$, and Galois extensions of $Z$ inside $D$ whose Galois groups are $p$-elementary abelian (Theorem \ref{EquivalenceGalTor}). In Section \ref{part2}, we use a specialization technique to establish our main result. For example, we prove that maximal subfields of $D(X)$ ``generically'' specialise to maximal subfields of $D$, and properties such as being Galois or purely inseparable over the centre also carry over generically (Theorem \ref{transfer}). \\

Applications are given in Section \ref{appli}. Let $L$ be any solvable Lie algebra, or a Zassenhaus algebra (see \ref{zassenhausdef} for the definition). It is proved that the enveloping skewfield of $L$ in characteristic $p > 2$ contains maximal subfields which are Galois (resp. purely inseparable of exponent 1) over the centre (Theorems \ref{env_galois} and \ref{zassenhausthm}). A crucial ingredient for the proof in the solvable case is the following. We prove that the enveloping field of a $p$-envelope of $L$ is isomorphic to a ring of rational functions over the enveloping skewfield of $L$ (Proposition \ref{reductiontorestricted}). In view of the previous results, this allows us to reduce to the case of restrictable Lie algebras, which is known by results of J. Schue \cite{SchuePaper}. \\

As a consequence of these theorems, we also show that the enveloping skewfield of $L$ defines an element of order $p$ in the Brauer group of its centre, when $L$ is solvable and non-abelian, or a Zassenhaus algebra. This suggests the following conjecture: \\

{\noindent \bf Conjecture.} Let $L$ be a non-abelian Lie algebra over a fields of characteristic $p > 2$, and let $K(L)$ be the enveloping skewfield. Then, $K(L)$ defines an element of order $p$ in the Brauer group of its centre. \\

{\noindent \bf Acknowledgements.} The authors wish to thank Rolf Farnsteiner for interesting comments and advice, and Amiram Braun for providing the original motivation for this work. They are also grateful to the referee for insightful suggestions.

\section{A reduction principle for division rings}

In what follows, we denote by $[V:D] := \dim_D(V)$ for a left vector space $V$ over a division ring $D$. For an algebra $A$, we denote $Z(A)$ the centre of $A$. For a prime number $p$, we denote by $\ZZ_p$ the cyclic group with $p$ elements and $\FF_p$ the field with $p$ elements; we use this notation to emphasise the field structure. \\

\subsection{Preliminaries: tori in $p$-division algebras} \label{part1}

\subsubsection{}
Before we deal with the reduction principle, we need some results on commutative subfields and tori in $p$-division algebras. Let $D$ be a $p$-division algebra, that is to say, a division ring of characteristic $p > 0$, of dimension some power of $p$ over its centre. We are interested in linking the notion of a  torus in $D$ with some class of subfields of $D$, which are Galois extensions of the centre $Z$. Recall that an element $t \in D$ is {\it toral} if $t^p - t \in Z$. Alternatively, this means that the inner derivation $\ad(t)$ is a toral element in the restricted Lie algebra $\Der_Z(D)$ \cite[p. 79]{SF}. A {\it torus} is a commutative $Z$-subspace $T \subseteq D$ which is spanned by toral elements. In particular, $\ad(T)$ is a torus in $\Der_Z(L)$ \cite[p.86]{SF}. We define the {\it rank of $T$} to be $[\ad(T) : Z]$.

Clearly, the unit element $1$ is toral, and if $T_0$ is a torus, then $Z + T_0$ is a torus as well, of same rank. Since we are concerned with the adjoint action of tori on the division ring $D$, we will henceforth only consider tori containing $1$.

\subsubsection{}
We recall some standard facts related to actions of a torus. Let $T$ be a torus, then there is a weight space decomposition
\begin{equation} \label{weightspace}
D = \bigoplus_{\lambda \in \Lambda} D_\lambda,
\end{equation}
where $\Lambda \subseteq T^*= \Hom_Z(T, Z)$ is the {\it set of weights (of $T$ in $D$)}. By definition,
\begin{equation*}
D_\lambda =  \left\{ x \in D\ |\ (\forall\, t \in T),\ [t,x] = \lambda(t) x \right\},
\end{equation*}
and $\Lambda$ is the set of linear forms $\lambda$ such that $D_\lambda \neq (0)$. Note that $D_0 = C_D(T)$, the centralizer of $T$ in $D$. It is easily seen that each $D_\lambda$ is a $D_0$-vector space (on the left and on the right), of dimension 1. Furthermore, one readily checks that $\Lambda$ is an additive subgroup of $T^*$, and the decomposition \eqref{weightspace} is a $\Lambda$-grading of $D$.

\subsubsection{} \label{GaloisTorus} The following result is essentially known \cite[Section 2]{SchuePaper}. We give a different proof and a more precise statement. \\

\begin{lemma}
Let $D$ be a $p$-division algebra with centre $Z$. Let $T\subseteq D$ be a torus of rank $d$, and $\Lambda$ be the corresponding set of weights. Then: 
\begin{enumerate}
\item The group $\Lambda \simeq \ZZ_p^d$.
\item Let $Z(T) \subseteq D$ be the subfield generated by $Z$ and $T$. Then $Z(T)$ is Galois over $Z$, and $\Gal (Z(T) / Z) \simeq \Lambda$. \\
\end{enumerate}
\end{lemma}

\proof
(1) We may assume that $T$ contains $1$. Let $\{t_0\ldots,t_d\}$ be a toral basis of $T$, with $t_0 = 1$. Let $T_p = \sum_{i=1}^d \FF_p t_i$, and $\Lambda_p:=\{ \lambda|_{T_p} \ |\ \lambda \in \LambdaÊ\} \subseteq \Hom_{\FF_p}(T_p,Z)$. Since $t_0 = 1$ acts trivially on $D$, it is clear that $\Lambda \simeq \Lambda_p$. We will show that $\Lambda_p = T_p^*:=\Hom_{\FF_p}(T_p,\FF_p)$, which will prove our first assertion.

First we show that $\Lambda_p \subseteq T_p^*$. For each $i \in \{Ê1,\ldots,d\}$, we have $(\ad\, t_i)^p - (\ad\, t_i) = 0$. Let $\lambda \in \Lambda$; since each $\lambda(t_i)$ is an eigenvalue of $\ad\, t_i$, we obtain $\lambda(t_i) \in \FF_p$ as we wanted. For the reverse inclusion, we consider the natural non-degenerate pairing
\begin{equation*}
\begin{array}{ccc}
T_p \times T_p^* & \to & \FF_p \\
(t,\lambda) & \mapsto & \lambda(t).
\end{array}
\end{equation*}
Let $t \in \Lambda_p^{\perp} \subseteq T_p$. For all $\lambda \in \Lambda_p$ and all $x_\lambda \in D_\lambda$, we have $[t,x_\lambda] = \lambda(t) x_\lambda = 0$. Owing to \eqref{weightspace}, we obtain $[t,D] = 0$, so that $t \in Z \cap T_p$. Since $\{1,t_1,\ldots,t_d\}$ is a $Z$-linearly independent family, we get $t=0$. This proves $\Lambda_p^{\perp} = (0)$, hence $\Lambda_p = T_p^*$. \\

(2) For all $i \in \{Ê1,\ldots,d\}$, we have $t_i^p-t_i \in Z$, hence $[Z(T):Z] \leq p^d$. Furthermore, since each $t_i$ is separable over $Z$, it follows that $Z(T)$ is separable over $Z$. In particular it admits a primitive element, say $\alpha \in Z(T)$. 

Let $P(X) \in Z[X]$ be the minimal polynomial of $\alpha$ over $Z$, so that $\deg(P) = [Z(T):Z]$. It is known that the cardinality $|\Aut_Z Z(T)|$ is the number of roots of $P(X)$ in $Z(T)$, whence $|\Aut_Z Z(T)| \leq [Z(T):Z]$. Thus, to show that $Z(T)$ is normal (and hence Galois) over $Z$ it suffices to prove that $|\Aut_Z Z(T)| = [Z(T):Z]$.

For all $\lambda \in \Lambda$, choose a non-zero element $x_\lambda \in D_\lambda$. For all $t \in T$, it is easily seen that $x_\lambda t  x_\lambda^{-1} = t - \lambda(t)$, so that the inner automorphism defined by $x_\lambda$ induces an automorphism $\sigma_\lambda \in \Aut_Z Z(T)$. One readily checks that the assignment $\lambda \in \Lambda \mapsto \sigma_\lambda \in \Aut_Z Z(T)$ is a group homomorphism. It is also injective, because $\sigma_\lambda = \id$ if and only if $ t - \lambda(t) = \sigma_\lambda(t) = t$ for all $t \in T$. It follows $p^d \geq [Z(T):Z] \geq |\Aut_Z Z(T)| \geq |\Lambda| = p^d$. Hence, equality holds everywhere. This shows that $Z(T)$ is Galois over $Z$, with $\Gal(Z(T)/Z) \simeq \Lambda$.

\subsubsection{\bf Theorem.} \label{EquivalenceGalTor} \it
Let $D$ be a finite-dimensional $p$-division algebra over its centre $Z$. Let $K \subseteq D$ be a commutative extension field of $Z$. The following are equivalent:
\begin{itemize}
\item[{\bf (i)}] There exists a torus $T \subseteq K$ of rank $d$, such that $K = Z(T)$;
\item[{\bf (ii)}] $K$ is a Galois extension of $Z$, and $\Gal(K/Z)$ is a $p$-elementary abelian group of rank $d$. \\
\end{itemize}

\rm 
\proof
{\bf (i) $\Rightarrow$ (ii)} follows from Lemma \ref{GaloisTorus}, as well as the equality of ranks. 

For {\bf (ii) $\Rightarrow$ (i)}, assume that $K$ is Galois over $Z$, with $\Gal(K/Z) \simeq \ZZ_p^d=:\Gamma$. For each $i \in \{Ê1,\ldots,d \}$, let $\Gamma_i = \ZZ_p \times \ldots \times \{Ê0 \} \times \ldots \times \ZZ_p$, where the trivial group occurs on the $i$-th slot. Set $K_i = K^{\Gamma_i}$. By \cite[Cor. VI.1.16]{Lang}, we have $K = K_1 \cdots K_d$. Furthermore, each $K_i$ is Galois over $Z$ with Galois group $\Gamma / \Gamma_i \simeq \ZZ_p$. By the Theorem of Artin-Schreier \cite[Thm. VI.6.4]{Lang}, there exist $c_i \in Z, t_i \in K_i$ such that $K_i = Z(t_i)$ and $t_i^p - t_i - c_i = 0$. Then, $T = \sum_{i = 1}^d Z t_i$ is a torus such that $K = Z(T)$.

\subsubsection{} \label{walrus}
We record a few more general results on tori in $D$. \\

\begin{prop}
Let $D$ be a division $p$-algebra with centre $Z$. Let $[D:Z] = p^{2n}$, and let $T$ be a torus of rank $d$. Recall the weight space decomposition $D = \bigoplus_{\lambda \in \Lambda} D_\lambda$. Then, the following are equivalent:
\begin{itemize}
\item[{\bf (i)}] $n = d$;
\item[{\bf (ii)}] $Z(T)$ is a maximal subfield;
\item[{\bf (iii)}] $D_0$ is a maximal subfield;
\item[{\bf (iv)}] $D_0 = Z(T)$;
\item[{\bf (v)}] $D_0$ is commutative;
\item[{\bf (vi)}] $|\Lambda| = p^n$. \\
\end{itemize}
\end{prop}

\proof
By Lemma \ref{GaloisTorus}, we have $|\Lambda| = p^d$, which proves {\bf (i) $\iff$ (vi)}. The same lemma also gives $[Z(T):Z]=p^d$. Recall that a commutative subfield $K \subseteq D$ containing $Z$ is maximal if and only if $[K:Z] = p^n$, yielding {\bf (i) $\iff$ (ii)}. Alternatively, such a field $K$ is maximal commutative if and only if $C_D(K) = K$, if and only if $C_D(K)$ is commutative. Taking into account the fact that $C_D(Z(T)) = C_D(T) = D_0$, we readily obtain {\bf (ii) $\iff$ (iv) $\iff$ (v)}. Finally, we have $C_D(D_0) \subseteq C_D(T) = D_0$. Hence, $C_D(D_0) = D_0$ if and only if $D_0$ is commutative, if and only if $C_D(D_0)$ is maximal commutative. This proves {\bf (iii) $\iff$ (v)}.

\subsection{The reduction principle} \label{part2}

\subsubsection{}

Let $Z$ be an infinite field of arbitrary characteristic. Consider $Z[X] = Z[u_1,\ldots,u_q]$ a polynomial ring in $q$ variables over $Z$, and $Z(X) = \Frac\, Z[X]$. Let $X = \Hom_Z(Z[X],X)$ be the set of algebra homomorphisms $Z[X] \to Z$. The evaluation maps provide an identification of $X$ with the affine space $\mathbb{A}^q(Z) = Z^q$. Namely, an element $(\lambda_1,\ldots,\lambda_q) \in \mathbb{A}^q(Z)$ is identified with the evaluation morphism 
\begin{equation*}
\lambda: \left\{ \begin{array}{rcl}
 Z[u_1,\ldots,u_q] & \to & Z \\
f(u_1,\ldots,u_q) & \mapsto & f(\lambda_1,\ldots,\lambda_q).
\end{array} \right.
\end{equation*}

A property dependent on a parametre $\lambda \in X$ is said to hold {\it generically}, or  {\it for almost all $\lambda \in X$}, if it holds on a non-empty open subset with respect to the Zariski topology of $X \equiv \mathbb{A}^q(Z)$. A typical example for an open subset is
\begin{equation*}
\Omega_f:= \{Ê\lambda \in X\, |\, \lambda(f) \neq 0 \},
\end{equation*}
for a given $f \in Z[X]$. The field $Z$ being infinite, we have $\Omega_f = \emptyset$ if and only if $f = 0$. Moreover, finite intersections of non-empty open subsets are non-empty. In particular, if finitely many properties hold generically, their conjunction also holds generically.

\subsubsection{Rational functions over a division ring} \label{rationalfunctions}
Let $D$ be a finite-dimensional central division algebra over $Z$. We consider $D[X]:=D \otimes_Z Z[X]$, a polynomial ring with coefficients in $D$. We will identify $D$ and $Z[X]$ with the subalgebras $D \otimes_Z Z$ and $Z \otimes_Z Z[X]$ of $D[X]$. The following results are well-known: \\

\begin{lemma}
\begin{enumerate}
\item The ring $D[X]$ has a skewfield of fractions, denoted $D(X)$. Further, $D(X) \simeq D \otimes_Z Z(X)$.
\item The centre of $D(X)$ is $Z(X)$, and $[D(X):Z(X)] = [D:Z]$. \\
\end{enumerate}
\end{lemma}

\subsubsection{} \label{independence}
Any element $\lambda \in X$ defines an algebra homomorphism $\pi_{\lambda} = \id \overline{\otimes} \lambda : D[X] \to D$. More precisely, for $d \otimes f \in D \otimes_Z Z[X]$, we have $\pi_{ \lambda} (d \otimes f) = d \otimes \lambda(f)$. For any subspace $V \subseteq D(X)$ over $Z(X)$, and any $\lambda \in X$, we define $\overline{V}_\lambda := \pi_\lambda(V \cap D[X]) \subseteq D$, the {\it specialization of $V$ along $\lambda$}. If $V \supseteq Z(X)$, then $\overline{V}_\lambda \supseteq Z$. We will need the following simple lemma: \\

\begin{lemma}
\begin{enumerate}
\item Let $Ê\{Êa_1,\ldots,a_n \} \subseteq D[X]$ be $Z(X)$-linearly independent. Then, there exists a dense open subset $\Omega \subseteq X$ such that the specializations $\{Ê\pi_\lambda(a_1),\ldots,\pi_\lambda(a_n)\} \subseteq D$ are linearly independent over $Z$, for all $\lambda \in \Omega$.
\item Let $V \subseteq D(X)$ be a $Z(X)$-subspace of dimension $n$. There exists a dense open subset $\Omega \subseteq X$ such that the specialization $\overline{V}_\lambda$ is a $Z$-subspace of dimension $n$, for all $\lambda \in \Omega$. \\
\end{enumerate}
\end{lemma}

\proof
(1) Let $\cB= \{\beta_1,\ldots,\beta_N\}$ be a basis of $D$ over $Z$, so that it is also a basis of $D[X]$ over $Z[X]$. Decompose each $a_i = \sum_{j = 1}^N f_{ij}\, \beta_j$, with $f_{ij} \in Z[X]$. Let $A = [f_{ij}] \in M_{n,N}\big( Z[X] \big)$. Since $\{ a_1,\ldots,a_n \}$ is linearly independent over $Z(X)$, there is an $n \times n$ submatrix $A'$ such that the minor $\det(A') \in Z[X] \smallsetminus\{0\}$.

Now note that for all $i$, $\pi_\lambda(a_i) = \sum_{j}  \lambda(f_{ij}) \beta_j $, so that the matrix $A_{\lambda} := [\lambda(f_{ij})]$ represents the vectors $\{Ê\pi_\lambda(a_1),\ldots,\pi_\lambda(a_n) \}$ in the basis $\cB$. Then, the $n \times n$ minor $\det \big((A')_{\lambda} \big)= \lambda(\det\, A')$, which is non-zero for almost all $\lambda \in X$. Hence, the matrix $A_{\lambda}$ has full rank for almost all $\lambda$, in which case the family $\{Ê\pi_\lambda(a_1),\ldots, \pi_\lambda(a_n)\}$ is linearly independent over $Z$. \\

(2) Let $\{Êa_1,\ldots,a_n\}$ be a $Z(X)$-basis of $V$. After multiplication by a suitable non-zero element of $Z[X]$, we may assume that all $a_i \in D[X]$. By (1), these elements almost always reduce to linearly independent elements in $\overline{V}_\lambda$, hence $[\overline{V}_\lambda : Z] \geq [V : Z(X)]$ for almost all $\lambda \in X$. Conversely, let $\{ b_1,\ldots,b_m \} \subseteq V \cap D[X]$ be a lift of some $Z$-basis of $\overline{V}_\lambda$ and let $B \in M_{m,N}\big( Z[X] \big)$ be the corresponding coefficients matrix. Since $\{Ê\pi_\lambda(b_1),\ldots, \pi_\lambda(b_m)\}$ are $Z$-linearly independent, as above there exists a non-vanishing $m \times m$ minor in the reduced matrix $B_{\lambda}$. The corresponding minor of $B$ is also nonzero, so that the matrix $B$ has rank $m$ over $Z(X)$. It readily follows $ [V:Z(X)] \geq [\overline{V}_\lambda:Z]$. 

\subsubsection{} \label{reduction}
Now we are ready to prove the reduction principle. We keep the previous notations. \\

\begin{thm} 
Let $K \subseteq D(X)$ be an extension field of $Z(X)$, and $\lambda \in X$.
\begin{enumerate}
\item The specialization $\overline{K}_\lambda \subseteq D$ is an extension field of $Z$, and for generic $\lambda \in Z$ we have $[\overline{K}_\lambda : Z] = [K : Z(X)]$. In particular, if $K$ is a maximal commutative subfield of $D(X)$, then $\overline{K}_\lambda$ is generically a maximal subfield of $D$.
\item If $K$ is Galois over $Z(X)$, then $K_{\lambda}$ is Galois over $Z$ for generic $\lambda \in X$.
\end{enumerate}

If $\Char(Z) = p > 0$, we have in addition:

\begin{enumerate}
\setcounter{enumi}{2}
\item If $K$ is purely inseparable of exponent $r$ over $Z(X)$, then $\overline{K}_\lambda$ is purely inseparable over $Z$, of exponent $\leq r$. Equality holds for generic $\lambda \in X$.
\item If $K$ is Galois over $Z(X)$, with Galois group $\ZZ_p^r$, then for generic $\lambda \in X$, $\overline{K}_\lambda$ is Galois over $Z$, with $\Gal(\overline{K}_\lambda / Z) \simeq \ZZ_p^r$. \\
\end{enumerate}
\end{thm}

\proof\ 
(1) By construction, $\overline{K}_\lambda$ is a finite-dimensional commutative domain over $Z$, so it is a field. Now, if $K \subseteq D(X)$ is a maximal subfield, then $[K:Z(X)]^2 = [D(X):Z(X)]$. By Lemma \ref{independence}, for almost all $\lambda \in X$, we have $[\overline{K}_\lambda:Z]^2 = [K:Z(X)]^2 = [D:Z]$, and hence $\overline{K}_\lambda$ is a maximal subfield of $D$. \\

(2) Choose an element $\alpha \in K$ which is primitive over $Z(X)$. After multiplying by a suitable element of $Z[X]$ we may assume that $\alpha \in D[X]$. Let $P(T) = \sum_{i = 0}^n c_i T^i \in Z(X)[T]$ be the minimal polynomial of $\alpha$ over $Z(X)$. Since $K$ is Galois over $Z(X)$, this polynomial splits into linear factors $P(T) = \prod_{i = 1}^n (T - \alpha_i)$, where each $\alpha_i \in K$. Now choose an element $c \in Z[X] \smallsetminus \{ 0 \}$ such that $c \alpha_i \in D[X]$ for all $i$. Then $\alpha$ is a root of $c^n P(T) = \prod_{i = 1}^n (cT - c \alpha_i) \in Z[X][T]$. For almost all $\lambda \in X$, the element $\lambda(c) \neq 0$. Then $\pi_\lambda(\alpha)$ is a root of
\begin{equation*}
P_\lambda(T) := \prod_{i = 1}^n \left( T - \lambda(c)^{-1} \pi_\lambda(c \alpha_i) \right) \in Z[T].
\end{equation*}
Indeed, since $\alpha \in D[X]$ we can write $\pi_\lambda(c \alpha) = \lambda(c) \pi_\lambda(\alpha)$, and hence $(T - \alpha) \mid P_\lambda(T)$.

Now note that $\{c \alpha_1,\ldots,c \alpha_n\}$ is a $Z(X)$-basis of $K$. By Lemma \ref{independence}, $\{ \pi_\lambda(c \alpha_1),\ldots,\pi_\lambda(c \alpha_n) \}$ is a $Z$-basis of $\overline{K}_\lambda$ for almost all $\lambda \in X$. In that case, $P_\lambda(T)$ is a separable polynomial, and $\overline{K}_\lambda = Z \big( \pi_\lambda(c \alpha_1),\ldots,\pi_\lambda(c \alpha_n) \big)$ is the splitting field of $P_\lambda(T)$. This proves that $\overline{K}_\lambda$ is a Galois extension of $Z$. \\

(3) Let $x \in \overline{K}_\lambda$, and choose an element $a \in K \cap D[X]$ with $\pi_\lambda(a) = x$. Since $K$ is purely inseparable over $Z(X)$, of exponent $r$, we have $a^{p^r} \in Z[X]$, hence, $x^{p^r} = \pi_\lambda(a^{p ^r}) \in Z$. 

We check that the inseparability exponents coincide for almost all $\lambda \in X$. There exists $a \in K$ such that $\{Ê1,a, a^p\ldots,a^{p ^{r-1}}\}$ is linearly independent over $Z(X)$. We may assume that $a \in D[X]$.  By Lemma \ref{independence}, for almost all $\lambda \in X$ the family $\{Ê1, \pi_\lambda(a), \ldots, \pi_\lambda(a)^{p^{r-1}}\}$ is linearly independent over $Z$: in that case, the inseparability exponent of $\overline{K}_\lambda$ over $Z$ is $> r-1$. \\

(4) Recall that being Galois with a $p$-elementary abelian Galois group is equivalent to being generated by toral elements (Theorem \ref{EquivalenceGalTor}). So we can write $K = Z(X)(t_1,\ldots,t_n)$, where the $t_i$ are toral and $\{Ê1,t_1,\ldots,t_n\}$ are $Z(X)$-linearly independent. It suffices to show that for almost all $\lambda \in X$, there exist toral elements $\tau_1,\ldots,\tau_n \in \overline{K}_\lambda$ such that $\{Ê1,\tau_1,\ldots,\tau_n\}$ are $Z$-linearly independent. Indeed, under these assumptions, we also know that $[K:Z(X)] = p^n = [Z(\tau_1,\ldots,\tau_n) : Z]$, yielding $ \overline{K}_\lambda = Z(\tau_1,\ldots,\tau_n)$.

There exists $c \in Z[X] \smallsetminus\{Ê0 \}$ such that all $c t_i \in D[X]$. For almost all $\lambda \in Z$, the family $\{ \pi_\lambda(c),Ê\pi_\lambda(ct_1), \ldots, \pi_\lambda(ct_n)\}$ is linearly independent over $Z$. In particular $\pi_\lambda(c) \neq 0$. A straightforward computation shows that each element $(c t_i)^p - c^{p-1} (c t_i)$ is central in $D[X]$. We obtain that each $\pi_\lambda(c t_i)^p - \pi_\lambda(c)^{p-1} \pi_\lambda(c t_i)$ is central, so that each $\tau_i := \pi_\lambda(c)^{-1} \pi_\lambda(c t_i)$ is toral in $\overline{K}_\lambda$. And by choice of $\lambda$, the family  $\{Ê1,\tau_1,\ldots,\tau_n\}$ is $Z$-linearly independent.

\subsubsection{Transfer theorems} \label{transfer}
Let $D$ be a finite-dimensional division algebra over its centre $Z$, and $D(X)$ be a division ring of rational functions in several variables over $D$. Recall that the {\it exponent} of a central simple algebra $R$, denoted by $\Exp(R)$, is the order of $R$ in the Brauer group of its centre $Z$ \cite[p. 214]{rowenBook}. In other words, this is the smallest integer $n \geq 1$ such that the $n$-th tensor power $R^{\otimes n}$ is isomorphic to some matrix algebra $M_N(Z)$ over $Z$. \\

\begin{thm}
\begin{enumerate}
\item One has $[D(X):Z(X)] = [D:Z]$, and $\Exp\, D(X) = \Exp\, D$.
\item $D(X)$ has a maximal subfield which is Galois over $Z(X)$ if and only if $D$ has a maximal subfield which is Galois over $Z$.
\item When $\Char(Z) = p > 0$: $D(X)$ has a maximal subfield which is purely inseparable of exponent $r$ (resp. Galois with Galois group $\ZZ_p^r$) over $Z(X)$ if and only if $D$ has a maximal subfield with the same property over $Z$. \\
\end{enumerate}
\end{thm}

\proof
(1) The identity $[D(X):Z(X)] = [D:Z]$ is stated in \ref{rationalfunctions}. For the exponent, recall that $D(X) \simeq D \otimes_Z Z(X)$ as algebras over $Z(X)$. For tensor powers, we compute:
\begin{eqnarray*}
D(X) \otimes_{Z(X)} D(X) & \simeq & D \otimes_Z Z(X) \otimes_{Z(X)} Z(X) \otimes_Z D \\
& \simeq & D \otimes_Z Z(X) \otimes_Z D \\
& \simeq & D \otimes_Z D \otimes_Z Z(X).
\end{eqnarray*}
We obtain inductively that $D(X)^{\otimes n} \simeq D^{\otimes n} \otimes_Z Z(X)$, where the tensor power on the left is taken over $Z(X)$ and the one on the right over $Z$. If $D^{\otimes n}$ is trivial in the Brauer group $\Br(Z)$, then $D(X)^{\otimes n}$ is trivial in $\Br \big( Z(X) \big)$, so $\Exp\, D(X) \mid \Exp(D)$. Conversely, assume that $D^{\otimes n} \otimes_Z Z(X) \simeq M_q\big( Z(X)\big)$, for some $q \geq 1$. We know that $D^{\otimes n} \simeq M_N(\Delta)$, for some central division $Z$-algebra $\Delta$ and some integer $N \geq 1$; so we have $M_N(\Delta) \otimes_Z Z(X) \simeq M_q \big( Z(X)\big)$. Using the fact that $M_N(\Delta) \otimes_Z Z(X) \simeq M_N \big( \Delta(X) \big)$, we obtain an isomorphism 
\begin{equation*}
M_N \left( \Delta(X) \right) \simeq M_q \left( Z(X) \right).
\end{equation*}
This implies that $\Delta(X) \simeq Z(X)$. In particular, $\Delta$ is commutative, whence $\Delta = Z$. Thus, the algebra $D^{\otimes n}$ is trivial in the Brauer group $\Br(Z)$. Therefore $\Exp (D) \mid \Exp\, D(X)$ as we wanted to show.

(2) and (3) The ``only if'' part follows from Theorem \ref{reduction}. Conversely, if $D$ has a maximal subfield $K$ satisfying any of the properties listed in (2) or (3), one checks that $K \otimes_Z Z(X) \subseteq D(X)$ is a maximal subfield with the same property.

\section{Applications} \label{appli}

\subsection{Enveloping skewfields of non-restricted Lie algebras}

\subsubsection{}
In this section, $\FF$ denotes an algebraically closed field of characteristic $p > 0$. Let $L$ be a finite dimensional Lie algebra over $\FF$. Let $U(L)$ be its enveloping algebra with centre $Z(L)$, and $K(L) =\Frac\, U(L)$ be the division ring of fractions of $U(L)$. We denote by $C(L)$ the centre of $K(L)$. The main result here is to show that when $L$ is solvable, there always exist maximal subfields of $K(L)$ which are Galois or purely inseparable of exponent one over $C(L)$. \\

\subsubsection{}
Recall that a Lie algebra is {\it restrictable} if there exists a map $x \in L \mapsto x^{[p]} \in L$, such that $(\ad\, x)^p = \ad(x^{[p]})$ for all $x \in L$. If $L$ is restrictable, one can choose this map with some additional properties which mimick the properties of associative $p$-th powers in an associative algebra. In that case, the map is called a {\it $p$-mapping}, and the pair $(L,[p])$ is called a {\it restricted Lie algebra}.  We don't write down explicitly these properties here, as they are quite technical and irrelevant in our situation; see \cite[Chap. 2]{SF} for a comprehensive account. \\

Finally, in the enveloping algebra of a restricted Lie algebra, the subalgebra 
\begin{equation*}
Z_p(L):=\FF \langle x^p - x^{[p]}\ |\ x \in L \rangle \subseteq U(L)
\end{equation*}
is contained in the centre of $U(L)$, and called the {\it $p$-centre of $U(L)$}.

\subsubsection{} \label{psequence}
We briefly recall the notion of a $p$-envelope, see \cite[Section 2.5]{SF} or \cite[Section 1.1]{Strade} for details. Let $L$ be embedded in a restricted Lie algebra $G$. The {\it $p$-envelope of $L$ in $G$}, denoted $L_{(p)}$, is the smallest restricted Lie subalgebra of $G$ containing $L$. Note that the structure of $L_{(p)}$ depends on the initial embedding. For example, if $L \subseteq \mathfrak{gl}(n)$, then the corresponding $p$-envelope $L_{(p)}$ is finite-dimensional. On the other hand, consider the natural inclusion $L \subseteq U(L)$, then the associated $p$-envelope is infinite-dimensional. \\

In the sequel, we will abuse terminology by referring to ``a $p$-envelope'' of a Lie algebra $L$. By this we will always mean the $p$-envelope of $L$ in some unspecified larger finite-dimensional restricted Lie algebra. Since $L$ is finite-dimensional, it always affords finite-dimensional $p$-envelopes \cite[Prop 2.5.3]{SF}. The following lemma provides a general description of how $L$ embeds into $L_{(p)}$: \\

\begin{lemma}
Let $L$ be a finite dimensional Lie algebra over $\FF$, and let $L_{(p)}$ be a finite-dimensional $p$-envelope of $L$. Then there exists a sequence of ideals $L_{(p)} = L_q \supseteq L_{q-1} \supseteq \cdots \supseteq L_0 = L$ such that, for all $i \in \{Ê1,\ldots,q \}$:
\begin{enumerate}
\item $L_i = \FF x_i + L_{i-1}$ for some $x_i \in L_i$,
\item there exists $y_{i-1} \in L_{i-1}$ with $y_{i-1}^{[p]} = x_i$. \\
\end{enumerate}
\end{lemma}

\proof
We construct the $L_i$ inductively. For $0 \leq i < q$, we have $L_i^{[p]} \not \subseteq L_i$, so there exists $y_i \in L_{i}$ such that $x_{i+1}:=y_i^{[p]} \not \in L_i$. Set $L_{i+1}:=L_i \oplus \FF x_{i+1}$. By construction these subspaces satisfy the desired conditions. And it is easy to see that they are ideals, cf. \cite[Prop. 2.1.3]{SF}.

\subsubsection{}
The following result describes the comparative ring structures of $U(L)$ and $U(L_{(p)})$. \\

\begin{prop} \label{reductiontorestricted}
Let $L$ be a finite dimensional Lie algebra over $\FF$, and $L_{(p)}$  be a $p$-envelope. Then $U(L_{(p)}) \simeq U(L)[u_1,\ldots,u_q]$, a polynomial extension in $q = \dim_\FF(L_{(p)} / L)$ variables. In particular, the enveloping field $K(L_{(p)})$ is isomorphic to a ring of rational functions over $K(L)$. \\
\end{prop}

\proof
Let $x_1,\ldots,x_q$, $y_0,\ldots,y_{q-1}$ and $L_0,\ldots,L_q$ be as in Lemma \ref{psequence}. Then, the elements $u_i:=x_i - y_{i-1}^p$ are in the $p$-centre of $U(L_{(p)})$, hence commute with $U(L)$. Now observe that each $u_i \equiv x_i \mod U(L_{i-1})$. It follows from the PBW theorem that the powers $\{1, u_i, u_i^2, \ldots \}$ form a basis of $U(L_i)$ as a module over $U(L_{i-1})$. Hence, the monomials in $u_1,\ldots,u_q$ are central and form a basis of $U(L_{(p)})$ as a $U(L)$-module, and the result readily follows.

\subsubsection{\bf Theorem.} \label{env_galois} \it
Assume that $\Char(\FF) = p > 2$. Let $L$ be a finite-dimensional non-abelian solvable Lie algebra over $\FF$. Then, the division ring $K(L)$ has the following properties:
\begin{enumerate}
\item There exists a maximal subfield $F \subseteq K(L)$ which is Galois over the centre and the Galois group is $p$-elementary abelian;
\item there exists a maximal subfield $E \subseteq K(L)$ which is purely inseparable, of exponent 1 over the centre;
\item $\Exp\, K(L) = p$. \\
\end{enumerate}

\rm
\proof
By Proposition \ref{reductiontorestricted} and Theorem \ref{transfer}, it is enough to show the properties for restricted solvable algebras. Then (1) follows from \cite[Thm. 3]{SchuePaper} and (2) follows from \cite [Thm. 2]{SchuePaper}. Property (3) follows from (2) and \cite[Thm. 4.1.8]{jacobson}.

\subsubsection{}
As a consequence of Theorem \ref{env_galois}, we obtain the following result. For a solvable Lie algebra in characteristic $p > 2$, there always exists a torus $T \subseteq K(L)$ which is ``maximal'' in the sense that $C_{K(L)}(T)$ is commutative. Alternatively, by Proposition  \ref{walrus}, this means $T$ has rank $n$, where $[K(L):C(L)] = p^{2n}$. If $L$ is restricted then it follows from Schue's results \cite{SchuePaper}. \\
%  However, in non restricted case, this is not trivial, since there might be no toral elements in $L$, as we can see in the following example. Let $L$ be 3-dimensional non restricted Lie algebra defined by $[x,y]=y$, $[x,z]=y+z$ and $[y,z]=0$. Then clearly no torus exists in $L$, but we can easily see that the maximal torus $S$ in $D(L)$ is spanned by $\FF_{p}x^{p} + \FF_{p}(x-x^{p})zy^{-1}$.

\begin{cor}
Assume that $\Char(\FF) = p > 2$. Let $L$ be a finite dimensional solvable Lie algebra over $\FF$. Then there exists a torus $T \subseteq K(L)$ such that $C_{K(L)}(T)$ is commutative.
\end{cor}

\subsection{The Zassenhaus algebra}

\subsubsection{} \label{zassenhausdef}
Let $\FF$ be algebraically closed of characteristic $p > 2$, and let $m \geq 1$ be a fixed integer. The {\it Zassenhaus algebra} is the simple Lie algebra of Cartan type $W(1,m)$ \cite[Chap. 4.2]{Strade}. Explicitly, $W(1,m)$ has a basis $\{ e_{-1}, e_0,\ldots,e_{p^m-2} \}$ with brackets:
\begin{equation*}
[e_i,e_j] = \left( \binom{i + j + 1}{i} - \binom{i + j + 1}{j} \right) e_{i+j},
\end{equation*}
so $[e_i,e_j] = 0$ when $i + j \not \in \{ -1,\ldots,p^m - 2 \}$.  \\

\subsubsection{\bf Theorem} \label{zassenhausthm} \it
Let $\FF$ be algebraically closed with $\Char(\FF) > 2$.  Then, the enveloping skewfield $K \big( W(1,m) \big)$ has the following properties:
\begin{enumerate}
\item There exists a maximal subfield which is Galois over the centre and the Galois group is $p$-elementary abelian;
\item there exists a maximal subfield which is purely inseparable, of exponent 1 over the centre;
\item $\Exp\, K \big( W(1,m) \big) = p$. \\
\end{enumerate}

\rm

\proof\
For ease of notation, let $L:=W(1,m)$. It is easy to see that the subspace $H:=\sum_{i \geq 0} \FF\, e_i$ is a solvable Lie subalgebra of codimension 1 in $L$. By Theorem \ref{env_galois}, the properties of the theorem are satisfied in $K(H)$. By \cite[Prop. 2]{ermolaev}, there exists a central element $z \in U(L)$ of the form $z = a e_{-1} + b$, where $a,b \in U(H)$, $a \neq 0$. Using the PBW theorem, we can see that $z$ is transcendental over $U(H)$. Furthermore, it is clear that $K(L) = \Frac\, U(H)[z] = K(H)(z)$, so applying the transfer Theorem \ref{transfer} yields the result.

\end{document}